%% file: Text-II.tex
\documentclass[11pt,twoside]{article}

\input{packages.tex}

\title{
  {\huge Complex Analysis of Real Functions \\[1.5ex]}
  II: Singular Schwartz Distributions }

\author{
  \Large Jorge L. deLyra\footnote{Email: delyra@latt.if.usp.br} \\
  Department of Mathematical Physics \\
  Physics Institute \\
  University of São Paulo }

\date{May 26, 2018}


\input{math-defs.tex}

\begin{document}\maketitle

\begin{abstract}
  \noindent
  In the context of the complex-analytic structure within the unit disk
  centered at the origin of the complex plane, that was presented in a
  previous paper, we show that singular Schwartz distributions can be
  represented within that same structure, so long as one defines the
  limits involved in an appropriate way. In that previous paper it was
  shown that essentially all integrable real functions can be represented
  within the complex-analytic structure. The infinite collection of
  singular objects which we analyze here can thus be represented side by
  side with those real functions, thus allowing all these objects to be
  treated in a unified way.
\end{abstract}

\section{Introduction}\label{Sec01}

In a previous paper~\cite{CAoRFI} we introduced a certain complex-analytic
structure within the unit disk of the complex plane, and showed that one
can represent essentially all integrable real functions within that
structure. In this paper we will show that one can represent within the
same structure the singular objects known as distributions, loosely in the
sense of the Schwartz theory of
distributions~\cite{DTSchwartz,DTLighthill}, which are also known as
generalized real functions. All these objects will be interpreted as parts
of this larger complex-analytic structure, within which they can be
treated and manipulated in a robust and unified way.

In Sections~\ref{Sec02} and~\ref{Sec03} we will establish the relation
between the complex-analytic structure and the singular distributions.
There we will show that one obtains these objects through the properties
of certain limits to the unit circle, involving a particular set of inner
analytic functions, which will be presented explicitly. Following what was
shown in~\cite{CAoRFI} for integrable real functions, each singular
distribution will be associated to a corresponding inner analytic
function. In fact, we will show that the entire space of all singular
Schwartz distributions defined within a compact domain is contained within
this complex-analytic structure. In Section~\ref{Sec04} we will analyze a
certain collection of integrable real functions which are closely related
to the singular distributions, through the concept of infinite
integral-differential chains of functions.

For ease of reference, we include here a one-page synopsis of the
complex-analytic structure introduced in~\cite{CAoRFI}. It consists of
certain elements within complex analysis~\cite{CVchurchill}, as well as of
their main properties.

\paragraph{Synopsis:} The Complex-Analytic Structure\\

\noindent
An {\em inner analytic function} $w(z)$ is simply a complex function which
is analytic within the open unit disk. An inner analytic function that has
the additional property that $w(0)=0$ is a {\em proper inner analytic
  function}. The {\em angular derivative} of an inner analytic function is
defined by

\noindent
\begin{equation}
  w^{\ldot}(z)
  =
  \ii
  z\,
  \frac{dw(z)}{dz}.
\end{equation}

\noindent
By construction we have that $w^{\ldot}(0)=0$, for all $w(z)$. The {\em
  angular primitive} of an inner analytic function is defined by

\begin{equation}
  w^{-1\ldot}(z)
  =
  -\ii
  \int_{0}^{z}dz'\,
  \frac{w(z')-w(0)}{z'}.
\end{equation}

\noindent
By construction we have that $w^{-1\ldot}(0)=0$, for all $w(z)$. In terms
of a system of polar coordinates $(\rho,\theta)$ on the complex plane,
these two analytic operations are equivalent to differentiation and
integration with respect to $\theta$, taken at constant $\rho$. These two
operations stay within the space of inner analytic functions, they also
stay within the space of proper inner analytic functions, and they are the
inverses of one another. Using these operations, and starting from any
proper inner analytic function $w^{0\ldot}(z)$, one constructs an infinite
{\em integral-differential chain} of proper inner analytic functions,

\begin{equation}
  \left\{
    \ldots,
    w^{-3\ldot}(z),
    w^{-2\ldot}(z),
    w^{-1\ldot}(z),
    w^{0\ldot}(z),
    w^{1\ldot}(z),
    w^{2\ldot}(z),
    w^{3\ldot}(z),
    \ldots\;
  \right\}.
\end{equation}

\noindent
Two different such integral-differential chains cannot ever intersect each
other. There is a {\em single} integral-differential chain of proper inner
analytic functions which is a constant chain, namely the null chain, in
which all members are the null function $w(z)\equiv 0$.

A general scheme for the classification of all possible singularities of
inner analytic functions is established. A singularity of an inner
analytic function $w(z)$ at a point $z_{1}$ on the unit circle is a {\em
  soft singularity} if the limit of $w(z)$ to that point exists and is
finite. Otherwise, it is a {\em hard singularity}. Angular integration
takes soft singularities to other soft singularities, and angular
differentiation takes hard singularities to other hard singularities.

Gradations of softness and hardness are then established. A hard
singularity that becomes a soft one by means of a single angular
integration is a {\em borderline hard} singularity, with degree of
hardness zero. The {\em degree of softness} of a soft singularity is the
number of angular differentiations that result in a borderline hard
singularity, and the {\em degree of hardness} of a hard singularity is the
number of angular integrations that result in a borderline hard
singularity. Singularities which are either soft or borderline hard are
integrable ones. Hard singularities which are not borderline hard are
non-integrable ones.

Given an integrable real function $f(\theta)$ on the unit circle, one can
construct from it a unique corresponding inner analytic function $w(z)$.
Real functions are obtained through the $\rho\to 1_{(-)}$ limit of the
real and imaginary parts of each such inner analytic function and, in
particular, the real function $f(\theta)$ is obtained from the real part
of $w(z)$ in this limit. The pair of real functions obtained from the real
and imaginary parts of one and the same inner analytic function are said
to be mutually Fourier-conjugate real functions.

Singularities of real functions can be classified in a way which is
analogous to the corresponding complex classification. Integrable real
functions are typically associated with inner analytic functions that have
singularities which are either soft or at most borderline hard. This ends
our synopsis.

\vspace{2.6ex}

\noindent
When we discuss real functions in this paper, some properties will be
globally assumed for these functions, just as was done in~\cite{CAoRFI}.
These are rather weak conditions to be imposed on these functions, that
will be in force throughout this paper. It is to be understood, without
any need for further comment, that these conditions are valid whenever
real functions appear in the arguments. These weak conditions certainly
hold for any integrable real functions that are obtained as restrictions
of corresponding inner analytic functions to the unit circle.

The most basic condition is that the real functions must be measurable in
the sense of Lebesgue, with the usual Lebesgue
measure~\cite{RARudin,RARoyden}. The second global condition we will
impose is that the functions have no removable singularities. The third
and last global condition is that the number of hard singularities on the
unit circle be finite, and hence that they be all isolated from one
another. There will be no limitation on the number of soft singularities.

The material contained in this paper is a development, reorganization and
extension of some of the material found, sometimes still in rather
rudimentary form, in the
papers~\cite{FTotCPI,FTotCPII,FTotCPIII,FTotCPIV,FTotCPV}.

\section{The Dirac Delta ``Function''}\label{Sec02}

This is where we begin the discussion of inner analytic functions that
have hard singularities with strictly positive degrees of hardness. Let us
start by simply introducing a certain particular inner analytic function
of $z$. If $z_{1}$ is a point on the unit circle, this function is defined
as a very simple rational function of $z$,

\begin{equation}\label{EQInnDelt}
  w_{\delta}(z,z_{1})
  =
  \frac{1}{2\pi}
  -
  \frac{1}{\pi}\,
  \frac{z}{z-z_{1}}.
\end{equation}

\noindent
This inner analytic function has a single point of singularity, which is a
simple pole at $z_{1}$. This is a hard singularity with degree of hardness
equal to one. Our objective here is to examine the properties of the real
part $u_{\delta}(\rho,\theta,\theta_{1})$ of this inner analytic function,

\begin{equation}
  w_{\delta}(z,z_{1})
  =
  u_{\delta}(\rho,\theta,\theta_{1})
  +
  \ii
  v_{\delta}(\rho,\theta,\theta_{1}).
\end{equation}

\noindent
We will prove that in the $\rho\to 1_{(-)}$ limit
$u_{\delta}(\rho,\theta,\theta_{1})$ can be interpreted as a {\em Schwartz
  distribution}~\cite{DTSchwartz,DTLighthill}, namely as the singular
object known as the {\em Dirac delta ``function''}, which we will denote
by $\delta(\theta-\theta_{1})$. This object is also known as a {\em
  generalized real function}, since it is not really a real function in
the usual sense of the term. In the Schwartz theory of distributions this
object plays the role of an integration kernel for a certain distribution.
Note that $w_{\delta}(z,z_{1})$ can, in fact, be written explicitly as a
function of $\rho$ and $\theta-\theta_{1}$. Since we have that
$z=\rho\exp(\ii\theta)$ and that $z_{1}=\exp(\ii\theta_{1})$, we have at
once that

\begin{equation}
  w_{\delta}(z,z_{1})
  =
  \frac{1}{2\pi}
  -
  \frac{1}{\pi}\,
  \frac
  {\rho\e{\iii(\theta-\theta_{1})}}
  {\rho\e{\iii(\theta-\theta_{1})}-1}.
\end{equation}

\noindent
The definition of the Dirac delta ``function'' is that it is a symbol for
a limiting process, which satisfies certain conditions. In our case here
the limiting process will be the limit $\rho\to 1_{(-)}$ from the open
unit disk to the unit circle. The limit of
$u_{\delta}(\rho,\theta,\theta_{1})$ represents the delta ``function'' in
the sense that it satisfies the conditions that follow.

\begin{enumerate}

\item The defining limit of $\delta(\theta-\theta_{1})$ tends to zero when
  one takes the $\rho\to 1_{(-)}$ limit while keeping
  $\theta\neq\theta_{1}$.

\item The defining limit of $\delta(\theta-\theta_{1})$ diverges to
  positive infinity when one takes the $\rho\to 1_{(-)}$ limit with
  $\theta=\theta_{1}$.

\item In the $\rho\to 1_{(-)}$ limit the integral

  \begin{equation}
    \int_{a}^{b}d\theta\,
    \delta(\theta-\theta_{1})
    =
    1
  \end{equation}

  \noindent
  has the value shown, for any open interval $(a,b)$ which contains the
  point $\theta_{1}$.

\item Given any continuous integrable function $g(\theta)$, in the
  $\rho\to 1_{(-)}$ limit the integral

  \begin{equation}
    \int_{a}^{b}d\theta\,
    g(\theta)
    \delta(\theta-\theta_{1})
    =
    g(\theta_{1})
  \end{equation}

  \noindent
  has the value shown, for any open interval $(a,b)$ which contains the
  point $\theta_{1}$.

  This is the usual form of this condition, when it is formulated in
  strictly real terms. However, we will impose a slight additional
  restriction on the real functions $g(\theta)$, by assuming that the
  limit to the point $z_{1}$ on the unit circle that corresponds to
  $\theta_{1}$, of the corresponding inner analytic function
  $w_{\gamma}(z)$, exists and is finite. This implies that $w_{\gamma}(z)$
  may have at $z_{1}$ a soft singularity, but not a hard singularity.

\end{enumerate}

\noindent
Note that, although it is customary to list both separately, the third
condition is in fact just a particular case of the fourth condition. It is
also arguable that the second condition is not really necessary, because
it is a consequence of the others. We may therefore consider that the only
really essential conditions are the first and the fourth ones.

The functions $g(\theta)$ are sometimes named {\em test functions} within
the Schwartz theory of distributions~\cite{DTSchwartz,DTLighthill}. The
additional part of the fourth condition, that the limit to the point
$z_{1}$ of the corresponding inner analytic function $w_{\gamma}(z)$ must
exist and be finite, consists of a weak limitation on these test
functions, and does not affect the definition of the singular distribution
itself. This is certainly the case for our definition here, since we
define this object through a definite and unique inner analytic function.

\vspace{2.6ex}

\noindent
In this section we will prove the following theorem.

\begin{theorem}\Colon\label{Theo01}
  The $\rho\to 1_{(-)}$ limit of the real part of the inner analytic
  function $w_{\delta}(z,z_{1})$ converges to the generalized function
  $\delta(\theta-\theta_{1})$.
\end{theorem}

\noindent
Before we attempt to prove this theorem, our first task is to write
explicitly the real and imaginary parts of $w_{\delta}(z,z_{1})$. In order
to do this we must now rationalize it,

\noindent
\begin{eqnarray}
  w_{\delta}(z,z_{1})
  & = &
  \frac{1}{2\pi}
  -
  \frac{1}{\pi}\,
  \frac{z(z^{*}-z_{1}^{*})}{(z-z_{1})(z^{*}-z_{1}^{*})}
  \nonumber\\
  & = &
  \frac{1}{2\pi}
  -
  \frac{1}{\pi}\,
  \frac
  {\rho^{2}-\rho\cos(\Delta\theta)-\ii\rho\sin(\Delta\theta)}
  {\rho^{2}-2\rho\cos(\Delta\theta)+1},
\end{eqnarray}

\noindent
where $\Delta\theta=\theta-\theta_{1}$. We must examine the real part of
this function,

\begin{equation}
  u_{\delta}(\rho,\theta,\theta_{1})
  =
  \frac{1}{2\pi}
  -
  \frac{1}{\pi}\,
  \frac
  {\rho\left[\rho-\cos(\Delta\theta)\right]}
  {\left(\rho^{2}+1\right)-2\rho\cos(\Delta\theta)}.
\end{equation}

\noindent
We are now ready to prove the theorem, which we will do by simply
verifying all the properties of the Dirac delta ``function''.

\begin{proof}\Colon
\end{proof}

\noindent
If we take the limit $\rho\to 1_{(-)}$, under the assumption that
$\Delta\theta\neq 0$, we get

\noindent
\begin{eqnarray}
  \lim_{\rho\to 1_{(-)}}u_{\delta}(\rho,\theta,\theta_{1})
  & = &
  \frac{1}{2\pi}
  -
  \frac{1}{\pi}\,
  \frac
  {1-\cos(\Delta\theta)}
  {2-2\cos(\Delta\theta)}
  \nonumber\\
  & = &
  0,
\end{eqnarray}

\noindent
which is the correct value for the case of the Dirac delta ``function''.
Thus we see that the first condition is satisfied.

\vspace{2.6ex}

\noindent
If, on the other hand, we calculate $u_{\delta}(\rho,\theta,\theta_{1})$
for $\Delta\theta=0$ and $\rho<1$ we obtain

\noindent
\begin{eqnarray}
  u_{\delta}(\rho,\theta_{1},\theta_{1})
  & = &
  \frac{1}{2\pi}
  -
  \frac{1}{\pi}\,
  \frac
  {\rho(\rho-1)}
  {(\rho-1)^{2}}
  \nonumber\\
  & = &
  \frac{1}{2\pi}
  -
  \frac{1}{\pi}\,
  \frac{\rho}{\rho-1},
\end{eqnarray}

\noindent
which diverges to positive infinity as $\rho\to 1_{(-)}$, as it should in
order to represent the singular Dirac delta ``function''. This establishes
that the second condition is satisfied.

\vspace{2.6ex}

\noindent
We now calculate the real integral of $u_{\delta}(\rho,\theta,\theta_{1})$
over the circle of radius $\rho<1$, which is given by

\noindent
\begin{eqnarray}\label{EQProp3}
  \int_{-\pi}^{\pi}d\theta\,
  \rho\,
  u_{\delta}(\rho,\theta,\theta_{1})
  & = &
  \frac{1}{2\pi}
  \int_{-\pi}^{\pi}d\theta\,
  \rho
  \left\{
    1
    -
    \frac
    {2\rho\left[\rho-\cos(\Delta\theta)\right]}
    {\left(\rho^{2}+1\right)-2\rho\cos(\Delta\theta)}
  \right\}
  \nonumber\\
  & = &
  \frac{\rho}{2\pi}
  \int_{-\pi}^{\pi}d(\Delta\theta)\,
  \frac
  {\left(1-\rho^{2}\right)}
  {\left(\rho^{2}+1\right)-2\rho\cos(\Delta\theta)}
  \nonumber\\
  & = &
  \frac{\left(1-\rho^{2}\right)}{4\pi}
  \int_{-\pi}^{\pi}d(\Delta\theta)\,
  \frac
  {1}
  {\left[\left(\rho^{2}+1\right)/(2\rho)\right]-\cos(\Delta\theta)},
\end{eqnarray}

\noindent
since $d(\Delta\theta)=d\theta$. This real integral over $\Delta\theta$
can be calculated by residues. We introduce an auxiliary complex variable
$\xi=\lambda\exp(\ii\Delta\theta)$, which becomes simply
$\exp(\ii\Delta\theta)$ on the unit circle $\lambda=1$. We have
$d\xi=\ii\xi d(\Delta\theta)$, and so we may write the integral on the
right-hand side as

\noindent
\begin{eqnarray}
  \int_{-\pi}^{\pi}d(\Delta\theta)\,
  \frac{1}{\left[\left(1+\rho^{2}\right)/(2\rho)\right]-\cos(\Delta\theta)}
  & = &
  \oint_{C}d\xi\,
  \frac{1}{\ii\xi}\,
  \frac{2}{\left[\left(1+\rho^{2}\right)/\rho\right]-\xi-1/\xi}
  \nonumber\\
  & = &
  2\ii
  \oint_{C}d\xi\,
  \frac{1}{1-\left[\left(1+\rho^{2}\right)/\rho\right]\xi+\xi^{2}},
\end{eqnarray}

\noindent
where the integral is now over the unit circle $C$ in the complex $\xi$
plane. The two roots of the quadratic polynomial on $\xi$ in the
denominator are given by

\noindent
\begin{eqnarray}
  \xi_{(+)}
  & = &
  1/\rho,
  \nonumber\\
  \xi_{(-)}
  & = &
  \rho.
\end{eqnarray}

\noindent
Since $\rho<1$, only the simple pole corresponding to $\xi_{(-)}$ lies
inside the integration contour, so we get for the integral

\noindent
\begin{eqnarray}
  \int_{-\pi}^{\pi}d\theta\,
  \frac{1}{\left[\left(1+\rho^{2}\right)/(2\rho)\right]-\cos(\Delta\theta)}
  & = &
  2\ii(2\pi\ii)
  \lim_{\xi\to\rho}
  \frac{1}{\xi-1/\rho}
  \nonumber\\
  & = &
  4\pi\,
  \frac{\rho}{\left(1-\rho^{2}\right)}.
\end{eqnarray}

\noindent
It follows that we have for the real integral in Equation~(\ref{EQProp3})

\noindent
\begin{eqnarray}
  \int_{-\pi}^{\pi}d\theta\,
  \rho\,
  u_{\delta}(\rho,\theta,\theta_{1})
  & = &
  \frac{\left(1-\rho^{2}\right)}{4\pi}\,
  4\pi\,
  \frac{\rho}{\left(1-\rho^{2}\right)}
  \nonumber\\
  & = &
  \rho,
\end{eqnarray}

\noindent
and thus we have that the integral is equal to $1$ in the $\rho\to
1_{(-)}$ limit. Once we have this result, and since according to the first
condition the integrand goes to zero everywhere on the unit circle except
at $\Delta\theta=0$, which is the same as $\theta=\theta_{1}$, the
integral can be changed to one over any open interval $(a,b)$ on the unit
circle containing the point $\theta_{1}$, without any change in its
limiting value. This establishes that the third condition is satisfied.

\vspace{2.6ex}

\noindent
In order to establish the validity of the fourth and last condition, we
consider an essentially arbitrary integrable real function $g(\theta)$,
with the additional restriction that it be continuous at the point
$z_{1}$. As was established in~\cite{CAoRFI}, it corresponds to an inner
analytic function

\begin{equation}
  w_{\gamma}(z)
  =
  u_{\gamma}(\rho,\theta)
  +
  \ii
  v_{\gamma}(\rho,\theta),
\end{equation}

\noindent
where we also assume that $g(\theta)$ is such that $w_{\gamma}(z)$ may
have at $z_{1}$ a soft singularity, but not a hard singularity, so that
its limit to $z_{1}$ exists. We now consider the following real
integral\footnote{Post-publication note: it is important to observe here
  that, from its very beginning, this argument goes through without any
  change if we replace $u_{\gamma}(\rho,\theta)$ directly by $g(\theta)$
  in the integrand. Hence, the theorem is true regardless of whether or
  not the real function $g(\theta)$ is represented by an inner analytic
  function.} over the circle of radius $\rho<1$,

\noindent
\begin{eqnarray}\label{EQProp4}
  \int_{-\pi}^{\pi}d\theta\,
  \rho\,
  u_{\gamma}(\rho,\theta)
  u_{\delta}(\rho,\theta,\theta_{1})
  & = &
  \frac{1}{2\pi}
  \int_{-\pi}^{\pi}d\theta\,
  \rho\,
  u_{\gamma}(\rho,\theta)
  \left\{
    1
    -
    \frac
    {2\rho\left[\rho-\cos(\Delta\theta)\right]}
    {\left(\rho^{2}+1\right)-2\rho\cos(\Delta\theta)}
  \right\}
  \nonumber\\
  & = &
  \Frac{\rho}{2\pi}
  \int_{-\pi}^{\pi}d(\Delta\theta)\,
  u_{\gamma}(\rho,\theta)\,
  \frac
  {\left(1-\rho^{2}\right)}
  {\left(\rho^{2}+1\right)-2\rho\cos(\Delta\theta)}
  \nonumber\\
  & = &
  \frac{\left(1-\rho^{2}\right)}{4\pi}
  \int_{-\pi}^{\pi}d(\Delta\theta)\,
  \frac
  {u_{\gamma}(\rho,\theta)}
  {\left[\left(\rho^{2}+1\right)/(2\rho)\right]-\cos(\Delta\theta)},
\end{eqnarray}

\noindent
since $d(\Delta\theta)=d\theta$. This real integral over $\Delta\theta$
can be calculated by residues, exactly like the one in
Equation~(\ref{EQProp3}) which appeared before in the case of the third
condition. The calculation is exactly the same except for the extra factor
of $u_{\gamma}(\rho,\theta)$ to be taken into consideration when
calculating the residue, so that we may write directly that

\noindent
\begin{eqnarray}
  \int_{-\pi}^{\pi}d(\Delta\theta)\,
  \frac
  {u_{\gamma}(\rho,\theta)}
  {\left[\left(\rho^{2}+1\right)/(2\rho)\right]-\cos(\Delta\theta)}
  & = &
  2\ii(2\pi\ii)
  \lim_{\xi\to\rho}
  \frac{u_{\gamma}(\rho,\theta)}{\xi-1/\rho}
  \nonumber\\
  & = &
  4\pi\,
  \frac{\rho}{\left(1-\rho^{2}\right)}
  \lim_{\xi\to\rho}u_{\gamma}(\rho,\theta).
\end{eqnarray}

\noindent
Note now that since $\xi=\lambda\exp(\ii\Delta\theta)$, and since we must
take the limit $\xi\to\rho$, we in fact have that in that limit

\begin{equation}
  \lambda\e{\iii\Delta\theta}
  =
  \rho,
\end{equation}

\noindent
which implies that $\lambda=\rho$ and that $\Delta\theta=0$, and therefore
that $\theta=\theta_{1}$. We must therefore write
$u_{\gamma}(\rho,\theta)$ at the point given by $\rho$ and $\theta_{1}$,
thus obtaining

\begin{equation}
  \int_{-\pi}^{\pi}d(\Delta\theta)\,
  \frac
  {u_{\gamma}(\rho,\theta)}
  {\left[\left(\rho^{2}+1\right)/(2\rho)\right]-\cos(\Delta\theta)}
  =
  4\pi\,
  \frac{\rho}{\left(1-\rho^{2}\right)}\,
  u_{\gamma}(\rho,\theta_{1}).
\end{equation}

\noindent
It follows that we have for the real integral in Equation~(\ref{EQProp4})

\noindent
\begin{eqnarray}
  \int_{-\pi}^{\pi}d\theta\,
  \rho\,
  u_{\gamma}(\rho,\theta)
  u_{\delta}(\rho,\theta,\theta_{1})
  & = &
  \frac{\left(1-\rho^{2}\right)}{4\pi}\,
  4\pi\,
  \frac{\rho}{\left(1-\rho^{2}\right)}\,
  u_{\gamma}(\rho,\theta_{1})
  \nonumber\\
  & = &
  \rho u_{\gamma}(\rho,\theta_{1}).
\end{eqnarray}

\noindent
Finally, we may now take the $\rho\to 1_{(-)}$ limit, since
$w_{\gamma}(z)$ and thus $u_{\gamma}(\rho,\theta)$ are well defined at
$z_{1}$ in that limit, and thus obtain

\noindent
\begin{eqnarray}
  \lim_{\rho\to 1_{(-)}}
  \int_{-\pi}^{\pi}d\theta\,
  \rho\,
  u_{\gamma}(\rho,\theta)
  u_{\delta}(\rho,\theta,\theta_{1})
  & = &
  u_{\gamma}(1,\theta_{1})
  \;\;\;\Rightarrow
  \nonumber\\
  \int_{-\pi}^{\pi}d\theta\,
  u_{\gamma}(1,\theta)
  \left[
    \lim_{\rho\to 1_{(-)}}
    u_{\delta}(\rho,\theta,\theta_{1})
  \right]
  & = &
  u_{\gamma}(1,\theta_{1})
  \;\;\;\Rightarrow
  \nonumber\\
  \int_{-\pi}^{\pi}d\theta\,
  g(\theta)
  \left[
    \lim_{\rho\to 1_{(-)}}
    u_{\delta}(\rho,\theta,\theta_{1})
  \right]
  & = &
  g(\theta_{1}),
\end{eqnarray}

\noindent
since $u_{\gamma}(\rho,\theta)$ converges to $g(\theta)$, in the $\rho\to
1_{(-)}$ limit, almost everywhere on the unit circle. Just as before, once
we have this result, and since according to the first condition the
integrand goes to zero everywhere on the unit circle except at
$\Delta\theta=0$, which is the same as $\theta=\theta_{1}$, the integral
can be changed to one over any open interval on the unit circle containing
the point $\theta_{1}$, without any change in its value. This establishes
that the fourth and last condition is satisfied.

\vspace{2.6ex}

\noindent
Having established all the properties, we may now write symbolically that

\begin{equation}
  \delta(\theta-\theta_{1})
  =
  \lim_{\rho\to 1_{(-)}}
  u_{\delta}(\rho,\theta,\theta_{1}).
\end{equation}

\noindent
This concludes the proof of Theorem~\ref{Theo01}.

\vspace{2.6ex}

\noindent
It is important to note that, when we adopt as the {\em definition} of the
Dirac delta ``function'' the $\rho\to 1_{(-)}$ limit of the real part of
the inner analytic function $w_{\delta}(z,z_{1})$, the limitations imposed
on the test functions $g(\theta)$ and on the corresponding inner analytic
functions $w_{\gamma}(z)$ become irrelevant. In fact, this definitions
stands by itself, and is independent of any set of test functions. Given
any integrable real function $f(\theta)$ and the corresponding inner
analytic function $w(z)$ with real part $u(\rho,\theta)$, we may always
assemble the real integral over a circle of radius $\rho<1$

\begin{equation}
  \int_{-\pi}^{\pi}d\theta\,
  \rho\,
  u(\rho,\theta)
  u_{\delta}(\rho,\theta,\theta_{1}),
\end{equation}

\noindent
which is {\em always} well defined within the open unit disk. It then
remains to be verified only whether or not the $\rho\to 1_{(-)}$ limit of
this integral exists, in order to define the corresponding integral

\begin{equation}
  \int_{-\pi}^{\pi}d\theta\,
  f(\theta)
  \delta(\theta-\theta_{1}).
\end{equation}

\noindent
This limit may exist for functions that do not satisfy the conditions
imposed on the test functions. In fact, one can do this for the real part
of {\em any} inner analytic function, regardless of whether or not it
corresponds to an integrable inner analytic function, so long as the
$\rho\to 1_{(-)}$ limit of $u(\rho,\theta)$ exists almost everywhere.
Whenever the $\rho\to 1_{(-)}$ limit of the integral exists, it defines
the action of the delta ``function'' on that particular real object.

It is also interesting to observe that the Dirac delta ``function'',
although it is not simply a conventional integrable real function, is in
effect an integrable real object, even if it corresponds to an inner
analytic functions that has a simple pole at $z_{1}$, which is a
non-integrable hard singularity, with degree of hardness equal to
one. This apparent contradiction is explained by the {\em orientation} of
the pole at $z=z_{1}$. If we consider the real part
$u_{\delta}(\rho,\theta)$ of the inner analytic function
$w_{\delta}(z,z_{1})$, although it is not integrable along curves arriving
at the singular point from most directions, there is one direction, that
of the unit circle, along which one can approach the singular point so
that $u_{\delta}(\rho,\theta)$ is identically zero during the approach,
which allows us to define its integral using the $\rho\to 1_{(-)}$
limit\footnote{Post-publication note: this characterizes the inner
  analytic function $w_{\delta}(z,z_{1})$ associated to the Dirac delta
  ``function'' as an {\em irregular} inner analytic function, since it is
  not integrable around its singular point, which is a hard singular point
  with degree of hardness $1$, while the corresponding real object is.}.
The same is {\em not} true, for example, for the imaginary part
$v_{\delta}(\rho,\theta)$ of the same inner analytic function, which
generates the Fourier-conjugate function to the delta ``function'', and
that diverges to infinity as $1/|z-z_{1}|$ when one approaches the
singular point along the unit circle, thus generating a non-integrable
real function in the $\rho\to 1_{(-)}$ limit.

In the development presented in~\cite{CAoRFI} the real functions were
represented by their Fourier coefficients, and the inner analytic
functions by their Taylor coefficients. We can easily do the same here, if
we observe that the inner analytic function $w_{\delta}(z,z_{1})$ in
Equation~(\ref{EQInnDelt}) is the sum of a geometric series,

\noindent
\begin{eqnarray}\label{EQSerDelt}
  w_{\delta}(z,z_{1})
  & = &
  \frac{1}{2\pi}
  +
  \frac{1}{\pi}\,
  \frac{z/z_{1}}{1-z/z_{1}}
  \nonumber\\
  & = &
  \frac{1}{2\pi}
  +
  \frac{1}{\pi}
  \sum_{k=1}^{\infty}
  \left(
    \frac{z}{z_{1}}
  \right)^{k}
  \nonumber\\
  & = &
  \frac{1}{2\pi}
  +
  \frac{1}{\pi}
  \sum_{k=1}^{\infty}
  \left[
    \cos(k\theta_{1})
    -
    \ii
    \sin(k\theta_{1})
  \right]
  z^{k}.
\end{eqnarray}

\noindent
This power series is the Taylor series of $w_{\delta}(z,z_{1})$ around the
origin, and therefore it follows that the Taylor coefficients of this
inner analytic function are given by

\noindent
\begin{eqnarray}
  c_{0}
  & = &
  \frac{1}{2\pi},
  \nonumber\\
  c_{k}
  & = &
  \frac{\cos(k\theta_{1})}{\pi}
  -
  \ii\,
  \frac{\sin(k\theta_{1})}{\pi},
\end{eqnarray}

\noindent
where $k\in\{1,2,3,\ldots,\infty\}$. Since according to the construction
presented in~\cite{CAoRFI} we have that $c_{0}=\alpha_{0}/2$ and that
$c_{k}=\alpha_{k}-\ii\beta_{k}$, we have for the Fourier coefficients of
the delta ``function''

\noindent
\begin{eqnarray}
  \alpha_{0}
  & = &
  \frac{1}{\pi},
  \nonumber\\
  \alpha_{k}
  & = &
  \frac{\cos(k\theta_{1})}{\pi},
  \nonumber\\
  \beta_{k}
  & = &
  \frac{\sin(k\theta_{1})}{\pi},
\end{eqnarray}

\noindent
where $k\in\{1,2,3,\ldots,\infty\}$. Note that these are in fact the
results one obtains via the integrals defining the Fourier
coefficients~\cite{FSchurchill},

\noindent
\begin{eqnarray}
  \alpha_{k}
  & = &
  \frac{1}{\pi}
  \int_{-\pi}^{\pi}d\theta\,
  \cos(k\theta)\delta(\theta-\theta_{1}),
  \nonumber\\
  \beta_{k}
  & = &
  \frac{1}{\pi}
  \int_{-\pi}^{\pi}d\theta\,
  \sin(k\theta)\delta(\theta-\theta_{1}),
\end{eqnarray}

\noindent
by simply using the fundamental property of the delta ``function''.

Having established the representation of the Dirac delta ``function''
within the structure of the inner analytic functions, in sequence we will
show that the Dirac delta ``function'' is not the only singular
distribution that can be represented by an inner analytic function. As we
will see, one can do the same for its first derivative, and in fact for
its derivatives of any order. This is an inevitable consequence of the
fact that the proper inner analytic function
$w_{\delta}^{0\ldot}(z,z_{1})$ associated to $w_{\delta}(z,z_{1})$ is a
member of an integral-differential chain.

\section{Derivatives of the Delta ``Function''}\label{Sec03}

The derivatives of the Dirac delta ``function'' are defined in a way which
is similar to that of the delta ``function'' itself. The first condition
is the same, and the second and third conditions are not really required.
The crucial difference is that the fourth condition in the definition of
the Dirac delta ``function'' is replaced by the second condition in the
list that follows. The ``function'' $\delta^{n\prime}(\theta-\theta_{1})$
is the $n^{\rm th}$ derivative of $\delta(\theta-\theta_{1})$ with respect
to $\theta$ if its defining limit $\rho\to 1_{(-)}$ satisfies the two
conditions that follow.

\begin{enumerate}

\item The defining limit of $\delta^{n\prime}(\theta-\theta_{1})$ tends to
  zero when one takes the $\rho\to 1_{(-)}$ limit while keeping
  $\theta\neq\theta_{1}$.

\item Given any integrable real function $g(\theta)$ which is
  differentiable to the $n^{\rm th}$ order, in the $\rho\to 1_{(-)}$ limit
  the integral

  \begin{equation}
    \int_{a}^{b}d\theta\,
    g(\theta)
    \delta^{n\prime}(\theta-\theta_{1})
    =
    (-1)^{n}g^{n\prime}(\theta_{1})
  \end{equation}

  \noindent
  has the value shown, for any open interval $(a,b)$ which contains the
  point $\theta_{1}$, where $g^{n\prime}(\theta)$ is the $n^{\rm th}$
  derivative of $g(\theta)$ with respect to $\theta$.

  This is the usual form of this condition, when it is formulated in
  strictly real terms. However, we will impose a slight additional
  restriction on the real functions $g(\theta)$, by assuming that the
  limit to the point $z_{1}$ on the unit circle that corresponds to
  $\theta_{1}$, of the $n^{\rm th}$ angular derivative of the
  corresponding inner analytic function $w_{\gamma}(z)$, exists and is
  finite. Since these proper inner analytic functions are all in the same
  integral-differential chain, this implies that the limits to $z_{1}$ of
  all the inner analytic functions $w^{m\ldot}_{\gamma}(z)$ exist, for all
  $0\leq m\leq n$.

\end{enumerate}

\noindent
The second condition above is, in fact, the fundamental property of each
derivative of the delta ``function'', including the ``function'' itself in
the case $n=0$. Just as in the case of the delta ``function'' itself, the
additional part of the second condition, involving the inner analytic
function $w_{\gamma}(z)$, consists of a weak limitation on the test
functions $g(\theta)$, and does not affect the definition of the singular
distributions themselves. This is certainly the case for our definitions
here, since we define each one of these objects through a definite and
unique inner analytic function.

\vspace{2.6ex}

\noindent
In this section we will prove the following theorem.

\begin{theorem}\Colon\label{Theo02}
  For every strictly positive integer $n$ there exists an inner analytic
  function $w_{\delta^{n\prime}}(z,z_{1})$ whose real part, in the
  $\rho\to 1_{(-)}$ limit, converges to
  $\delta^{n\prime}(\theta-\theta_{1})$.
\end{theorem}

\noindent
Before we attempt to prove this theorem, let us note that the proof relies
on a property of angular differentiation, which was established
in~\cite{CAoRFI}, namely that angular differentiation is equivalent to
partial differentiation with respect to $\theta$, at constant $\rho$. When
we take the $\rho\to 1_{(-)}$ limit, this translates to the fact that
taking the angular derivative of the inner analytic function $w(z)$ within
the open unit disk corresponds to taking the derivative with respect to
$\theta$, on the unit circle, of the corresponding real object.

If this derivative cannot be taken directly on the unit circle, then one
can {\em define} it by taking the angular derivative of the corresponding
inner analytic function and then considering the $\rho\to 1_{(-)}$ limit
of the real part of the resulting function. Since analytic functions can
be differentiated any number of times, the procedure can then be iterated
in order to define all the higher-order derivatives with respect to
$\theta$ on the unit circle. Equivalently, one can just consider traveling
along the integral-differential chain indefinitely in the differentiation
direction.

Consider therefore the integral-differential chain of proper inner
analytic functions that is obtained from the proper inner analytic
function associated to $w_{\delta}(z,z_{1})$, that is, the unique
integral-differential chain to which the proper inner analytic function

\noindent
\begin{eqnarray}\label{EQPropDelt}
  w_{\delta}^{0\ldot}(z,z_{1})
  & = &
  w_{\delta}(z,z_{1})
  -
  \frac{1}{2\pi}
  \nonumber\\
  & = &
  -\,
  \frac{1}{\pi}\,
  \frac{z}{z-z_{1}}
\end{eqnarray}

\noindent
belongs. Consider in particular the set of proper inner analytic functions
which is obtained from $w_{\delta}^{0\ldot}(z,z_{1})$ in the
differentiation direction along this chain, for which we have

\noindent
\begin{eqnarray}
  w_{\delta}^{n\ldot}(z,z_{1})
  & = &
  u_{\delta}^{n\prime}(\rho,\theta,\theta_{1})
  +
  \ii
  v_{\delta}^{n\prime}(\rho,\theta,\theta_{1})
  \nonumber\\
  & = &
  \frac{\partial^{n}}{\partial\theta^{n}}
  u_{\delta}(\rho,\theta,\theta_{1})
  +
  \ii\,
  \frac{\partial^{n}}{\partial\theta^{n}}
  v_{\delta}(\rho,\theta,\theta_{1}),
\end{eqnarray}

\noindent
for all strictly positive $n$, where we recall that

\begin{equation}
  w_{\delta}(z,z_{1})
  =
  u_{\delta}(\rho,\theta,\theta_{1})
  +
  \ii
  v_{\delta}(\rho,\theta,\theta_{1}).
\end{equation}

\noindent
We will now prove that in the $\rho\to 1_{(-)}$ limit we have

\begin{equation}\label{EQDerDel}
  \delta^{n\prime}(\theta-\theta_{1})
  =
  \lim_{\rho\to 1_{(-)}}
  u_{\delta}^{n\prime}(\rho,\theta,\theta_{1}),
\end{equation}

\noindent
for $n\in\{1,2,3,\ldots,\infty\}$, or, equivalently, that we have for the
inner analytic function $w_{\delta^{n\prime}}(z,z_{1})$ associated to the
derivative $\delta^{n\prime}(\theta-\theta_{1})$

\begin{equation}
  w_{\delta^{n\prime}}(z,z_{1})
  =
  w_{\delta}^{n\ldot}(z,z_{1}),
\end{equation}

\noindent
for $n\in\{1,2,3,\ldots,\infty\}$. We are now ready to prove the theorem,
as stated in Equation~(\ref{EQDerDel}). Let us first prove, however, that
the first condition holds for all the derivatives of the delta
``function''.

\newpage

\begin{proof}\Colon
\end{proof}

\noindent
Since $w_{\delta}(z,z_{1})$ has a {\em single} singular point at $z_{1}$,
the same is true for all its angular derivatives. Therefore the $\rho\to
1_{(-)}$ limit of all the angular derivatives exists everywhere within the
open interval of the unit circle that excludes the point
$\theta_{1}$. Since $u_{\delta}(1,\theta,\theta_{1})$ is identically zero
within this interval, and since angular differentiation within the open
unit disk corresponds to differentiation with respect to $\theta$ on the
unit circle, so that we have

\begin{equation}
  u_{\delta}^{n\prime}(1,\theta,\theta_{1})
  =
  \frac{\partial^{n}}{\partial\theta^{n}}
  u_{\delta}(1,\theta,\theta_{1}),
\end{equation}

\noindent
for $n\in\{1,2,3,\ldots,\infty\}$, it follows at once that

\noindent
\begin{eqnarray}
  u_{\delta}^{n\prime}(1,\theta,\theta_{1})
  & = &
  0
  \;\;\;\Rightarrow
  \nonumber\\
  \lim_{\rho\to 1_{(-)}}
  u_{\delta}^{n\prime}(\rho,\theta,\theta_{1})
  & = &
  0,
\end{eqnarray}

\noindent
for $n\in\{1,2,3,\ldots,\infty\}$, everywhere but at the singular point
$\theta_{1}$, for all values of $n$. This establishes that the first
condition holds.

\vspace{2.6ex}

\noindent
Let us now prove that the second condition, which relates directly to the
singular point, holds, leading to the result as stated in
Equation~(\ref{EQDerDel}).

\begin{proof}\Colon
\end{proof}

\noindent
In order to do this, we start with the case $n=1$, and consider the
following real integral on the circle of radius $\rho<1$, which we
integrate by parts, noting that the integrated term is zero because we are
integrating on a circle,

\begin{equation}
  \int_{-\pi}^{\pi}d\theta\,
  u_{\gamma}(\rho,\theta)
  \left[
    \frac{\partial}{\partial\theta}
    u_{\delta}(\rho,\theta,\theta_{1})
  \right]
  =
  -
  \int_{-\pi}^{\pi}d\theta\,
  \left[
    \frac{\partial}{\partial\theta}
    u_{\gamma}(\rho,\theta)
  \right]
  u_{\delta}(\rho,\theta,\theta_{1}),
\end{equation}

\noindent
where $w_{\gamma}(z)=u_{\gamma}(\rho,\theta)+\ii v_{\gamma}(\rho,\theta)$
is the inner analytic function associated to $g(\theta)$. Note that the
partial derivatives involved certainly exist, since both
$u_{\delta}(\rho,\theta,\theta_{1})$ and $u_{\gamma}(\rho,\theta)$ are the
real parts of inner analytic functions. If we now take the $\rho\to
1_{(-)}$ limit, on the right-hand side we recover the Dirac delta
``function'' on the unit circle, and therefore we have

\noindent
\begin{eqnarray}
  \int_{-\pi}^{\pi}d\theta\,
  g(\theta)
  \left[
    \lim_{\rho\to 1_{(-)}}
    \frac{\partial}{\partial\theta}
    u_{\delta}(\rho,\theta,\theta_{1})
  \right]
  & = &
  -
  \int_{-\pi}^{\pi}d\theta\,
  \left[
    \frac{d}{d\theta}
    g(\theta)
  \right]
  \delta(\theta-\theta_{1})
  \nonumber\\
  & = &
  (-1)\,
  g'(\theta_{1}),
\end{eqnarray}

\noindent
so long as $g(\theta)$ is differentiable, were we used the fundamental
property of the Dirac delta ``function''. We thus obtain the relation for
the derivative of the delta ``function'',

\begin{equation}
  \int_{-\pi}^{\pi}d\theta\,
  g(\theta)
  \delta'(\theta-\theta_{1})
  =
  (-1)\,
  g'(\theta_{1}),
\end{equation}

\noindent
where

\begin{equation}
  \delta'(\theta-\theta_{1})
  =
  \lim_{\rho\to 1_{(-)}}
  \frac{\partial}{\partial\theta}
  u_{\delta}(\rho,\theta,\theta_{1}).
\end{equation}

\noindent
We may therefore write that

\begin{equation}
  \delta'(\theta-\theta_{1})
  =
  \lim_{\rho\to 1_{(-)}}
  u'_{\delta}(\rho,\theta,\theta_{1}).
\end{equation}

\noindent
In this way we have obtained the result for $\delta'(\theta-\theta_{1})$
by using the known result for $\delta(\theta-\theta_{1})$. We may now
repeat this procedure to obtain the result for
$\delta^{2\prime}(\theta-\theta_{1})$ from the result for
$\delta'(\theta-\theta_{1})$, and therefore from the result for
$\delta(\theta-\theta_{1})$. We simply consider the following real
integral on the circle of radius $\rho<1$, which we again integrate by
parts, recalling that the integrated term is zero because we are
integrating on a circle,

\begin{equation}
  \int_{-\pi}^{\pi}d\theta\,
  u_{\gamma}(\rho,\theta)
  \left[
    \frac{\partial}{\partial\theta}
    u'_{\delta}(\rho,\theta,\theta_{1})
  \right]
  =
  -
  \int_{-\pi}^{\pi}d\theta\,
  \left[
    \frac{\partial}{\partial\theta}
    u_{\gamma}(\rho,\theta)
  \right]
  u'_{\delta}(\rho,\theta,\theta_{1}).
\end{equation}

\noindent
If we now take the $\rho\to 1_{(-)}$ limit, on the right-hand side we
recover the first derivative of the Dirac delta ``function'' on the unit
circle, and therefore we have

\noindent
\begin{eqnarray}
  \int_{-\pi}^{\pi}d\theta\,
  g(\theta)
  \left[
    \lim_{\rho\to 1_{(-)}}
    \frac{\partial}{\partial\theta}
    u'_{\delta}(\rho,\theta,\theta_{1})
  \right]
  & = &
  -
  \int_{-\pi}^{\pi}d\theta\,
  \left[
    \frac{d}{d\theta}
    g(\theta)
  \right]
  \delta'(\theta-\theta_{1})
  \nonumber\\
  & = &
  (-1)^{2}g^{2\prime}(\theta_{1}),
\end{eqnarray}

\noindent
so long as $g(\theta)$ is differentiable to second order, were we used the
fundamental property of the first derivative of the Dirac delta
``function''. We thus obtain the relation for the second derivative of the
delta ``function'',

\begin{equation}
  \int_{-\pi}^{\pi}d\theta\,
  g(\theta)
  \delta^{2\prime}(\theta-\theta_{1})
  =
  (-1)^{2}g^{2\prime}(\theta_{1}),
\end{equation}

\noindent
where

\begin{equation}
  \delta^{2\prime}(\theta-\theta_{1})
  =
  \lim_{\rho\to 1_{(-)}}
  \frac{\partial^{2}}{\partial\theta^{2}}
  u_{\delta}(\rho,\theta,\theta_{1}).
\end{equation}

\noindent
We may therefore write that

\begin{equation}
  \delta^{2\prime}(\theta-\theta_{1})
  =
  \lim_{\rho\to 1_{(-)}}
  u_{\delta}^{2\prime}(\rho,\theta,\theta_{1}).
\end{equation}

\noindent
Clearly, this procedure can be iterated $n$ times, thus resulting in the
relation

\begin{equation}
  \delta^{n\prime}(\theta-\theta_{1})
  =
  \lim_{\rho\to 1_{(-)}}
  u_{\delta}^{n\prime}(\rho,\theta,\theta_{1}),
\end{equation}

\noindent
for $n\in\{1,2,3,\ldots,\infty\}$. Note that all the derivatives with
respect to $\theta$ involved in the argument exist, for arbitrarily high
orders, since both $u_{\delta}(\rho,\theta,\theta_{1})$ and
$u_{\gamma}(\rho,\theta)$ are the real parts of inner analytic functions,
and thus are infinitely differentiable on both arguments.

We may now formalize the proof using finite induction. We thus assume the
results for the case $n-1$,

\noindent
\begin{eqnarray}
  \delta^{(n-1)\prime}(\theta-\theta_{1})
  & = &
  \lim_{\rho\to 1_{(-)}}
  u_{\delta}^{(n-1)\prime}(\rho,\theta,\theta_{1}),
  \nonumber\\
  \int_{a}^{b}d\theta\,
  g(\theta)
  \delta^{(n-1)\prime}(\theta-\theta_{1})
  & = &
  (-1)^{n-1}g^{(n-1)\prime}(\theta_{1}),
\end{eqnarray}

\noindent
and proceed to examine the next case. We consider therefore the following
real integral on the circle of radius $\rho<1$, which we integrate by
parts, recalling once more that the integrated term is zero because we are
integrating on a circle,

\begin{equation}
  \int_{-\pi}^{\pi}d\theta\,
  u_{\gamma}(\rho,\theta)
  \left[
    \frac{\partial}{\partial\theta}
    u_{\delta}^{(n-1)\prime}(\rho,\theta,\theta_{1})
  \right]
  =
  -
  \int_{-\pi}^{\pi}d\theta\,
  \left[
    \frac{\partial}{\partial\theta}
    u_{\gamma}(\rho,\theta)
  \right]
  u_{\delta}^{(n-1)\prime}(\rho,\theta,\theta_{1}).
\end{equation}

\noindent
If we now take the $\rho\to 1_{(-)}$ limit, on the right-hand side we
recover the $(n-1)^{\rm th}$ derivative of the Dirac delta ``function'' on
the unit circle, and therefore we have

\noindent
\begin{eqnarray}
  \int_{-\pi}^{\pi}d\theta\,
  g(\theta)
  \left[
    \lim_{\rho\to 1_{(-)}}
    \frac{\partial}{\partial\theta}
    u_{\delta}^{(n-1)\prime}(\rho,\theta,\theta_{1})
  \right]
  & = &
  -
  \int_{-\pi}^{\pi}d\theta\,
  \left[
    \frac{d}{d\theta}
    g(\theta)
  \right]
  \delta^{(n-1)\prime}(\theta-\theta_{1})
  \nonumber\\
  & = &
  (-1)^{n}\,
  g^{n\prime}(\theta_{1}),
\end{eqnarray}

\noindent
so long as $g(\theta)$ is differentiable to order $n$, were we used the
fundamental property of the $(n-1)^{\rm th}$ derivative of the Dirac delta
``function''. We thus obtain the relation for the $n^{\rm th}$ derivative
of the delta ``function'',

\begin{equation}
  \int_{-\pi}^{\pi}d\theta\,
  g(\theta)
  \delta^{n\prime}(\theta-\theta_{1})
  =
  (-1)^{n}\,
  g^{n\prime}(\theta_{1}),
\end{equation}

\noindent
where

\begin{equation}
  \delta^{n\prime}(\theta-\theta_{1})
  =
  \lim_{\rho\to 1_{(-)}}
  \frac{\partial^{n}}{\partial\theta^{n}}
  u_{\delta}(\rho,\theta,\theta_{1}).
\end{equation}

\noindent
We may therefore write that, by finite induction,

\begin{equation}
  \delta^{n\prime}(\theta-\theta_{1})
  =
  \lim_{\rho\to 1_{(-)}}
  u_{\delta}^{n\prime}(\rho,\theta,\theta_{1}),
\end{equation}

\noindent
for $n\in\{1,2,3,\ldots,\infty\}$. We have therefore completed the proof
of Theorem~\ref{Theo02}.

\vspace{2.6ex}

\noindent
It is important to note that, just as in the case of the Dirac delta
``function'', when we adopt as the {\em definition} of the $n^{\rm th}$
derivative of the delta ``function'' the $\rho\to 1_{(-)}$ limit of the
real part of the inner analytic function $w_{\delta}^{n\ldot}(z,z_{1})$,
for $n\in\{1,2,3,\ldots,\infty\}$, the limitations imposed on the test
functions $g(\theta)$ and on the corresponding inner analytic functions
$w_{\gamma}(z)$ become irrelevant. In fact, these definitions stand by
themselves, and are independent of any set of test functions. Not only one
can use them for any inner analytic functions derived from integrable real
functions, but one can do this for {\em any} inner analytic function
$w(z)$, regardless of whether or not it corresponds to an integrable real
function, so long as the $\rho\to 1_{(-)}$ limit of the corresponding real
part $u(\rho,\theta)$ exists almost everywhere. Just as in the case of the
Dirac delta ``function'', whenever the $\rho\to 1_{(-)}$ limit of the real
integral

\begin{equation}
  \int_{-\pi}^{\pi}d\theta\,
  \rho\,
  u(\rho,\theta)
  u_{\delta}^{n\prime}(\rho,\theta,\theta_{1}),
\end{equation}

\noindent
exists, it defines the action of the $n^{\rm th}$ derivative of the delta
``function'' on that particular real object.

It is also interesting to observe that, just as in the case of the Dirac
delta ``function'', it is true that its derivatives of all orders,
although they are not simply integrable real functions, are in fact
integrable real objects, even if they are related to inner analytic
functions with non-integrable hard singularities. Just as is the case for
the inner analytic function associated to the delta ``function'' itself,
the poles of the proper inner analytic functions associated to the
derivatives are always oriented in such a way that one can approach the
singularities along the unit circle while keeping the real parts of the
functions equal to zero, a fact that allows one to define the integrals on
$\theta$ of the real parts via the $\rho\to 1_{(-)}$
limit\footnote{Post-publication note: this characterizes all the inner
  analytic functions $w_{\delta}^{n\ldot}(z,z_{1})$ associated to the
  derivatives of the Dirac delta ``function'' as {\em irregular} inner
  analytic functions, since they are not integrable around their singular
  points, which are hard singular points with degrees of hardness $n+1$,
  while the corresponding real objects are.}. Just as in the case of the
delta ``function'', the Fourier-conjugate functions of the derivatives are
simply non-integrable real functions. This fact provides the first hint
that there must be some relation of such non-integrable real functions
with corresponding inner analytic functions.

In the development presented in~\cite{CAoRFI} the real functions were
represented by their Fourier coefficients, and the inner analytic
functions by their Taylor coefficients. The same can be done in our case
here. Starting from the power series for $w_{\delta}^{0\ldot}(z,z_{1})$
given in Equation~(\ref{EQSerDelt}), we can see that the definition of the
angular derivative implies that we have for the inner analytic functions
associated to the derivatives of the delta ``function'',

\begin{equation}
  w_{\delta}^{n\ldot}(z,z_{1})
  =
  \frac{1}{\pi}
  \sum_{k=1}^{\infty}
  \ii^{n}
  k^{n}
  \left[
    \cos(k\theta_{1})
    -
    \ii
    \sin(k\theta_{1})
  \right]
  z^{k},
\end{equation}

\noindent
for $n\in\{1,2,3,\ldots,\infty\}$, so that the corresponding Taylor
coefficients are given by $c_{0}^{(n)}=0$ and

\begin{equation}
  c_{k}^{(n)}
  =
  \frac{\ii^{n}k^{n}}{\pi}
  \left[
    \cos(k\theta_{1})
    -
    \ii\,
    \sin(k\theta_{1})
  \right],
\end{equation}

\noindent
for $n\in\{1,2,3,\ldots,\infty\}$, and where
$k\in\{1,2,3,\ldots,\infty\}$. The identification of the Fourier
coefficients $\alpha_{k}^{(n)}$ and $\beta_{k}^{(n)}$ will now depend on
the parity of $n$.

Once we have the Dirac delta ``function'' and all its derivatives, both as
inner analytic functions and as the corresponding real objects, we may
consider collections of such singular objects, with their singularities
located at all the possible points of the periodic interval $[-\pi,\pi]$,
as well as arbitrary linear combinations of some or all of them. There is
a well-known theorem of the Schwartz theory of
distributions~\cite{DTSchwartz,DTLighthill} which states that any
distribution which is singular at a given point $\theta_{1}$ can be
expressed as a linear combination of the Dirac delta ``function''
$\delta(\theta-\theta_{1})$ and its derivatives of arbitrarily high orders
$\delta^{n\prime}(\theta-\theta_{1})$.

Since, as was observed in~\cite{CAoRFI}, the set of all inner analytic
functions forms a vector space over the field of complex numbers, it is
immediately apparent that we may assemble such linear combinations within
the space of inner analytic functions. Therefore, the set of distributions
formed by the delta ``functions'' and all their derivatives, as defined
here, with their singularities located at all possible points of the unit
circle, constitutes a complete basis that spans the space of all possible
singular Schwartz distributions defined in a compact domain. We may
conclude therefore that the whole space of Schwartz distributions in a
compact domain is contained within the set of inner analytic functions.

It is interesting to note that, since we have the inner analytic function
that corresponds to the delta ``function'' in explicit form, we are in a
position to perform simple calculations in order to obtain in explicit
form the inner analytic functions that correspond to the first few
derivatives of the delta ``function''. For example, a few simple and
straightforward calculations lead to the following proper inner analytic
functions,

\noindent
\begin{eqnarray}
  w_{\delta^{1\prime}}(z,z_{1})
  & = &
  -\,
  \frac{1}{\pi\ii^{1}}\,
  \frac
  {zz_{1}}
  {(z-z_{1})^{2}},
  \nonumber\\
  w_{\delta^{2\prime}}(z,z_{1})
  & = &
  -\,
  \frac{1}{\pi\ii^{2}}\,
  \frac
  {z(z+z_{1})z_{1}}
  {(z-z_{1})^{3}},
  \nonumber\\
  w_{\delta^{3\prime}}(z,z_{1})
  & = &
  -\,
  \frac{1}{\pi\ii^{3}}\,
  \frac
  {z\left(z^{2}+4zz_{1}+z_{1}^{2}\right)z_{1}}
  {(z-z_{1})^{4}}.
\end{eqnarray}

\noindent
These proper inner analytic functions are all very simple rational
functions of the complex variable $z$, which can be written as functions
of only $z/z_{1}$, and hence as functions of only $\rho$ and
$\theta-\theta_{1}$. Note that we can induce from these examples that the
$n^{\rm th}$ derivative of the delta ``function'' is indeed represented by
an inner analytic function with a pole of order $n+1$ on the unit circle,
which is thus a hard singularity with degree of hardness $n+1$, as one
would expect from the structure of the corresponding integral-differential
chain.

\section{Piecewise Polynomial Functions}\label{Sec04}

It is important to note that the Dirac delta ``function'' and all its
derivatives, with their singularities located at a given point $z_{1}$ on
the unit circle, are all contained within a single integral-differential
chain, making up, in fact, only a part of that chain, the semi-infinite
chain starting from the delta ``function'' and propagating indefinitely in
the differentiation direction along the chain. However, the chain
propagates to infinity in both directions. In order to complete its
analysis, we must now determine what is the character of the real objects
in the remaining part of that chain, in the integration direction. In
fact, they are just integrable real functions, although they do have a
specific character. They consist of sections of polynomials wrapped around
the unit circle, of progressively higher orders, and progressively
smoother across the singular point, as functions of $\theta$, as one goes
along the integral-differential chain in the integration direction.
 
Let us illustrate this fact with a few simple examples. Instead of
performing angular integrations of the inner analytic functions, we will
do this by performing integrations on the unit circle. As was established
in~\cite{CAoRFI}, one can determine these functions by simple integration
on $\theta$, so long as one remembers two things: first, to make sure that
the real functions or related objects to be integrated on $\theta$ have
zero average over the unit circle, and second, to choose the integration
constant so that the resulting real functions also have zero average over
the unit circle. For example, the integral of the zero-average delta
``function''

\begin{equation}
  \delta^{0\prime}(\theta-\theta_{1})
  =
  \delta(\theta-\theta_{1})
  -
  \frac{1}{2\pi},
\end{equation}

\noindent
which is obtained from the real part of the proper inner analytic function
in Equation~(\ref{EQPropDelt}), can be integrated by means of the simple
use of the fundamental property of the delta ``function'', thus yielding

\noindent
\begin{equation}
  %
  \renewcommand{\arraystretch}{2.0}
  \begin{array}{rclcc}
    \delta^{-1\prime}(\Delta\theta)
    & = &
    \FFrac{1}{2}
    -
    \FFrac{\Delta\theta}{2\pi}
    &
    \mbox{for}
    &
    \Delta\theta>0,
    \nonumber\\
    \delta^{-1\prime}(\Delta\theta)
    & = &
    -\,
    \FFrac{1}{2}
    -
    \FFrac{\Delta\theta}{2\pi}
    &
    \mbox{for}
    &
    \Delta\theta<0,
  \end{array}
\end{equation}

\noindent
where $\Delta\theta=\theta-\theta_{1}$. This is a sectionally linear
function, with a single section consisting of the intervals $[-\pi,0)$ and
$(0,\pi]$, thus excluding the point $\Delta\theta=0$ where the singularity
lies, and with a unit-height step discontinuity at that point. Note that
it is an odd function of $\Delta\theta$. The next case can now be
calculated by straightforward integration, which yields

\begin{equation}
  %
  \renewcommand{\arraystretch}{2.0}
  \begin{array}{rclcc}
    \delta^{-2\prime}(\Delta\theta)
    & = &
    -\,
    \FFrac{\pi}{6}
    +
    \FFrac{\Delta\theta}{2}
    -
    \FFrac{\Delta\theta^{2}}{4\pi}
    &
    \mbox{for}
    &
    \Delta\theta>0,
    \nonumber\\
    \delta^{-2\prime}(\Delta\theta)
    & = &
    -\,
    \FFrac{\pi}{6}
    -
    \FFrac{\Delta\theta}{2}
    -
    \FFrac{\Delta\theta^{2}}{4\pi}
    &
    \mbox{for}
    &
    \Delta\theta<0.
  \end{array}
\end{equation}

\noindent
This is a sectionally quadratic function, this time a continuous function,
again with the same single section excluding the point $\Delta\theta=0$,
but now with a point of non-differentiability there. Note that it is an
even function of $\Delta\theta$. The next case yields, once more by
straightforward integration,

\noindent
\begin{equation}
  %
  \renewcommand{\arraystretch}{2.0}
  \begin{array}{rclcc}
    \delta^{-3\prime}(\Delta\theta)
    & = &
    -\,
    \FFrac{\pi\Delta\theta}{6}
    +
    \FFrac{\Delta\theta^{2}}{4}
    -
    \FFrac{\Delta\theta^{3}}{12\pi}
    &
    \mbox{for}
    &
    \Delta\theta>0,
    \nonumber\\
    \delta^{-3\prime}(\Delta\theta)
    & = &
    -\,
    \FFrac{\pi\Delta\theta}{6}
    -
    \FFrac{\Delta\theta^{2}}{4}
    -
    \FFrac{\Delta\theta^{3}}{12\pi}
    &
    \mbox{for}
    &
    \Delta\theta<0.
  \end{array}
\end{equation}

\noindent
This is a sectionally cubic continuous and differentiable function, again
with the same single section excluding the point $\Delta\theta=0$. Note
that it is an odd function of $\Delta\theta$. The trend is now quite
clear. All the real functions in the chain, in the integration direction
starting from the delta ``function'', are what we may call {\em piecewise
  polynomials}, even if we have just a single piece within a single
section of the unit circle, as is the case here. The $n^{\rm th}$ integral
is a piecewise polynomial of order $n$, which has zero average over the
unit circle, and which becomes progressively smoother across the singular
point as one goes along the integral-differential chain in the integration
direction.

In order to generalize this analysis, we must now consider linear
superpositions of delta ``functions'' and derivatives of delta
``functions'', with their singularities situated at various points on the
unit circle. A simple example of such a superposition, which we may use to
illustrate what happens when we make one, is that of two delta
``functions'', with singularities at $\theta=0$ and at $\theta=\pm\pi$,
added together with opposite signs,

\begin{equation}
  f(\theta)
  =
  \delta(\theta)
  -
  \delta(\theta-\pi),
\end{equation}

\noindent
that corresponds to the following inner analytic function, which this time
is already a proper inner analytic function, with two simple poles at
$z=\pm 1$,

\noindent
\begin{eqnarray}
  w(z)
  & = &
  -\,
  \frac{1}{\pi}\,
  \frac{z}{z-1}
  +
  \frac{1}{\pi}\,
  \frac{z}{z+1}
  \nonumber\\
  & = &
  -\,
  \frac{2}{\pi}\,
  \frac{z}{z^{2}-1}.
\end{eqnarray}

\noindent
Since we have now two singular points, one at $z=1$ and another at $z=-1$,
corresponding respectively to $\theta=0$ and $\theta=\pm\pi$, we have now
two sections, one in $(-\pi,0)$ and another in $(0,\pi)$. The inner
analytic functions at the integration side of the integral-differential
chain to which this function belongs are obtained by simply adding the
corresponding inner analytic functions at the integration sides of the
integral-differential chains of the two functions that are superposed. The
same is true for the corresponding real objects within each section of the
unit circle. Since the real functions corresponding to each one of the two
delta ``functions'' that were superposed are zero-average piecewise
polynomials, so are the real functions corresponding to the superposition.
For example, it is not difficult to show that the first integral is the
familiar square wave, with amplitude $1/2$,

\noindent
\begin{equation}
  %
  \renewcommand{\arraystretch}{2.0}
  \begin{array}{rclcc}
    f^{-1\prime}(\theta)
    & = &
    \FFrac{1}{2}
    &
    \mbox{for}
    &
    \theta>0,
    \nonumber\\
    f^{-1\prime}(\theta)
    & = &
    -\,
    \FFrac{1}{2}
    &
    \mbox{for}
    &
    \theta<0,
  \end{array}
\end{equation}

\noindent
which is a piecewise linear function with two sections, having unit-height
step discontinuities with opposite signs at the two singular points
$\theta=0$ and $\theta=\pm\pi$.

We want to determine what is the character of the real functions, in the
integration side of the resulting integral-differential chain, in the most
general case, when we consider arbitrary linear superpositions of a finite
number of delta ``functions'' and derivatives of delta ``functions'', with
their singularities situated at various points on the unit circle. From
the examples we see that, when we superpose several singular distributions
with their singularities at various points, the complete set of all the
singular points defines a new set of sections. Given one of these singular
points, since at least one of the distributions being superposed is
singular at that point, in general so is the superposition. Let there be
$N\geq 1$ singular points $\{\theta_{1},\ldots,\theta_{N}\}$ in the
superposition. It follows that in general we end up with a set of $N$
contiguous sections, consisting of open intervals between singular points,
that can be represented as the sequence

\begin{equation}
  \left\{
    \rule{0em}{2ex}
    (\theta_{1},\theta_{2}),
    \ldots,
    (\theta_{i-1},\theta_{i}),
    (\theta_{i},\theta_{i+1}),
    \ldots,
    (\theta_{N},\theta_{1})
  \right\},
\end{equation}

\noindent
where we see that the sequence goes around the unit circle, and where we
adopt the convention that each section $(\theta_{i},\theta_{i+1})$ is
numbered after the singular point $\theta_{i}$ at its left end. In
addition to this, since for each one of the distributions being superposed
the real functions on the integration side of the integral-differential
chain of the corresponding delta ``function'' are piecewise polynomials,
and since the sum of any finite number of polynomials is also a
polynomial, so are the real functions of the integral-differential chain
to which the superposition belongs, if we are at a point in that
integral-differential chain where all singular distributions have already
been integrated out. Let us establish a general notation for these
piecewise polynomial real functions, as well as a formal definition for
them.

\begin{definition}\Colon\label{Def01}
  Piecewise Polynomial Real Functions
\end{definition}

\noindent
Given a real function $f_{(n)}(\theta)$ that is defined in a piecewise
fashion by polynomials in $N\geq 1$ sections of the unit circle, with the
exclusion of a finite set of $N$ singular points $\theta_{i}$, with
$i\in\{1,\ldots,N\}$, so that the polynomial $P_{i}^{(n_{i})}(\theta)$ at
the $i^{\rm th}$ section has order $n_{i}$, we denote the function by

\begin{equation}
  f_{(n)}(\theta)
  =
  \left\{
    P_{i}^{(n_{i})}(\theta),
    i\in\{1,\ldots,N\}
  \right\},
\end{equation}

\noindent
where $n$ is the largest order among all the $N$ orders $n_{i}$. We say
that $f_{(n)}(\theta)$ is a {\em piecewise polynomial real function} of
order $n$.

\vspace{2.6ex}

\noindent
Note that, being made out of finite sections of polynomials, the real
function $f_{(n)}(\theta)$ is always an integrable real function. In fact,
it is also analytic within each section, so that the $N$ singularities
described above are the {\em only} singularities involved. Since
$f_{(n)}(\theta)$ is an integrable real function, let $w(z)$ be the inner
analytic function that corresponds to this integrable real function, as
constructed in~\cite{CAoRFI}. The $(n+1)^{\rm th}$ angular derivative of
$w(z)$ is the inner analytic function $w^{(n+1)\ldot}(z)$, which
corresponds therefore to the $(n+1)^{\rm th}$ derivative of
$f_{(n)}(\theta)$ with respect to $\theta$, that we denote by
$f_{(n)}^{(n+1)\prime}(\theta)$.

\vspace{2.6ex}

\noindent
In this section we will prove the following theorem.

\begin{theorem}\Colon\label{Theo03}
  If the real function $f_{(n)}(\theta)$ is a non-zero piecewise
  polynomial function of order $n$, defined in $N\geq 1$ sections of the
  unit circle, with the exclusion of a finite non-empty set of $N$
  singular points $\theta_{i}$, then and only then the derivative
  $f_{(n)}^{(n+1)\prime}(\theta)$ is the superposition of a non-empty set
  of delta ``functions'' and derivatives of delta ``functions'' on the
  unit circle, with the singularities located at some of the points
  $\theta_{i}$, and of nothing else.
\end{theorem}

\newpage

\begin{proof}\Colon
\end{proof}

\noindent
In order to prove this, first let us note that the derivative
$f_{(n)}^{(n+1)\prime}(\theta)$ is identically zero within all the open
intervals defining the sections. This is so because the maximum order of
all the polynomials involved is $n$, and the $(n+1)^{\rm th}$ derivative
of a polynomial of order equal to or less than $n$ is identically zero,

\begin{equation}
  f_{(n)}^{(n+1)\prime}(\theta)
  =
  0
  \mbox{\ \ \ for all\ \ \ }
  \theta\neq\theta_{i},
  i\in\{1,\ldots,N\}.
\end{equation}

\noindent
We conclude, therefore, that the real object represented by the inner
analytic function $w^{(n+1)\ldot}(z)$ has support only at the $N$ isolated
singular points $\theta_{i}$, thus implying that it can contain only
singular distributions.

Second, let us prove that the derivative cannot be identically zero over
the whole unit circle. In order to do this we note that one cannot have a
non-zero piecewise polynomial real function of order $n$, such as the one
described above, that is also continuous and differentiable to the order
$n$ on the whole unit circle. This is so because this hypothesis would
lead to an impossible integral-differential chain.

If this were possible, then starting from a non-zero real function
$f_{(n)}(\theta)$ that corresponds to a non-zero inner analytic function
$w(z)$, and after a finite number $n$ of steps along the differentiation
direction of the corresponding integral-differential chain, one would
arrive at a real function that is {\em continuous} over the whole unit
circle, that is {\em constant within each section} and that has {\em zero
  average} over the whole unit circle. It follows that such a function
would have to be identically zero, thus corresponding to the inner
analytic function $w(z)\equiv 0$. But this is not possible, because this
inner analytic function belongs to another chain, the one which is
constant, all members being $w(z)\equiv 0$, and we have shown
in~\cite{CAoRFI} that two different integral-differential chains cannot
intersect.

It follows that $f_{(n)}(\theta)$ can be globally differentiable at most
to order $n-1$, so that the $n^{\rm th}$ derivative is a discontinuous
function, and therefore the $(n+1)^{\rm th}$ derivative already gives rise
to singular distributions. Therefore, every real function that is
piecewise polynomial on the unit circle, of order $n$, when differentiated
$n+1$ times, so that it becomes zero within the open intervals
corresponding to the existing sections, will always result in the
superposition of some non-empty set of singular distributions with their
singularities located at the points between two consecutive sections.

We can also establish that {\em only} functions of this form give rise to
such superpositions of singular distributions and of nothing else. The
necessity of the fact that the real functions on integral-differential
chains generated by superpositions of singular distributions must be
piecewise polynomials comes directly from the fact that all such
distributions and all such superpositions of distributions are zero almost
everywhere, in fact everywhere but at their singular points. Due to this,
it is necessary that these real functions, upon a {\em finite} number
$n+1$ of differentiations, become zero everywhere strictly within the
sections, that is, within the open intervals between two successive
singularities. Therefore, within each open interval the condition over the
sectional real function at that interval is that

\begin{equation}
  \frac{d^{(n+1)}}{d\theta^{(n+1)}}f_{i}(\theta)
  \equiv
  0,
\end{equation}

\noindent
and the general solution of this ordinary differential equation of order
$n+1$ is a polynomial of order $n_{i}\leq n$, containing at most $n+1$
non-zero arbitrary constants,

\begin{equation}
  f_{i}(\theta)
  =
  P_{i}^{(n_{i})}(\theta).
\end{equation}

\noindent
Since only polynomials have the property of becoming identically zero
after a finite number of differentiations, it is therefore an absolute
necessity that these real functions be polynomials within each one of the
sections. This completes the proof of Theorem~\ref{Theo03}.

\vspace{2.6ex}

\noindent
Note that the inner analytic function $w^{(n+1)\ldot}(z)$ corresponding to
$f_{(n)}^{(n+1)\prime}(\theta)$ represents therefore the superposition of
a non-empty set of singular distributions with their singularities located
at the singular points. In other words, after $n+1$ angular
differentiations of $w(z)$, which correspond to $n+1$ straight
differentiations with respect to $\theta$ of the polynomials
$P_{i}^{(n_{i})}(\theta)$ within the sections, one is left with an inner
analytic function $w^{(n+1)\ldot}(z)$ whose real part converges to zero in
the $\rho\to 1_{(-)}$ limit, at all points on the unit circle which are
not one of the $N$ singular points.

It is interesting to note that, since we have the inner analytic function
that corresponds to the Dirac delta ``function'' in explicit form, it is
not difficult to calculate directly its first angular primitive. A few
simple and straightforward calculations lead to

\noindent
\begin{eqnarray}
  w_{\delta}^{-1\ldot}(z,z_{1})
  & = &
  \frac{\ii}{\pi}
  \int_{0}^{z}dz'\,
  \frac{z'}{z'-z_{1}}
  \nonumber\\
  & = &
  \frac{\ii}{\pi}
  \ln\!\left(\frac{z_{1}-z}{z_{1}}\right).
\end{eqnarray}

\noindent
This inner analytic function has a logarithmic singularity at $z_{1}$,
which is a borderline hard singularity. Note that, as expected, we have
that $w_{\delta}^{-1\ldot}(0,z_{1})=0$.

\section{Products of Distributions}\label{Sec05}

In the Schwartz theory of distributions one important theorem states that
it is not possible to define, in a general way, the product of two
distributions~\cite{NMSchwartz}, which has the effect that the space of
Schwartz distributions cannot be promoted from a vector space to an
algebra. In this section we will interpret this important fact in the
context of the representation of integrable real functions and singular
distributions in terms of inner analytic functions. We start by noting
that, although it is always possible to define the product of two inner
analytic functions, which is always an inner analytic function itself,
this does {\em not} correspond to the product of the two corresponding
real functions or related objects. If we have two inner analytic functions
given by

\noindent
\begin{eqnarray}\label{EQTwoInn}
  w_{1}(z)
  & = &
  u_{1}(\rho,\theta)
  +
  \ii
  v_{1}(\rho,\theta),
  \nonumber\\
  w_{2}(z)
  & = &
  u_{2}(\rho,\theta)
  +
  \ii
  v_{2}(\rho,\theta),
\end{eqnarray}

\noindent
the product of the two inner analytic functions is given by

\noindent
\begin{eqnarray}
  w(z)
  & = &
  \left[
    u_{1}(\rho,\theta)
    u_{2}(\rho,\theta)
    -
    v_{1}(\rho,\theta)
    v_{2}(\rho,\theta)
  \right]
  \nonumber\\
  &   &
  +
  \ii
  \left[
    u_{1}(\rho,\theta)
    v_{2}(\rho,\theta)
    +
    v_{1}(\rho,\theta)
    u_{2}(\rho,\theta)
  \right],
\end{eqnarray}

\noindent
whose real part is {\em not} just the product
$u_{1}(\rho,\theta)u_{2}(\rho,\theta)$. In fact, the problem of finding an
inner analytic function whose real part is this quantity often has no
solution. One can see this very simply by observing that both
$u_{1}(\rho,\theta)$ and $u_{2}(\rho,\theta)$ are always harmonic
functions, and that the product of two harmonic functions in general is
{\em not} a harmonic function. Since the real and imaginary parts of an
inner analytic function are always harmonic functions, it follows that the
problem posed in this way cannot be solved in general. The only simple
case in which we can see that the problem has a solution is that in which
one of the two functions being multiplied is a constant function.

Let us state in a general way the problem of the definition of the product
of two distributions. Suppose that we have two inner analytic functions
such as those in Equation~(\ref{EQTwoInn}). The two corresponding real
objects are $u_{1}(1,\theta)$ and $u_{2}(1,\theta)$, and their product,
assuming that it can be defined in strictly real terms, is simply the real
object $u_{1}(1,\theta)u_{2}(1,\theta)$. The problem of finding an inner
analytic function that corresponds to this product is the problem of
finding an harmonic function $u_{\pi}(\rho,\theta)$ whose limit to the
unit circle results in this real object,

\begin{equation}
  u_{\pi}(1,\theta)
  =
  u_{1}(1,\theta)u_{2}(1,\theta).
\end{equation}

\noindent
If one can find such a harmonic function, then it is always possible to
find its harmonic conjugate function $v_{\pi}(\rho,\theta)$ and therefore
to determine the inner analytic function

\begin{equation}
  w_{\pi}(z)
  =
  u_{\pi}(\rho,\theta)
  +
  \ii
  v_{\pi}(\rho,\theta),
\end{equation}

\noindent
which corresponds to the product of the two real objects. According to the
construction presented in~\cite{CAoRFI}, this can always be done so long
as the product $u_{1}(1,\theta)u_{2}(1,\theta)$ is an integrable real
function of $\theta$. However, if $u_{1}(1,\theta)$ and $u_{2}(1,\theta)$
are singular objects that can only be defined as limits from within the
open unit disk, then the product may not be definable in strictly real
terms, and it may not be possible to find an inner analytic function such
that the $\rho\to 1_{(-)}$ limit of its real part results in this product,
interpreted in some consistent way. This is the content of the theorem
that states that this cannot be done in general.

It is not too difficult to give examples of products which are not well
defined. It suffices to consider the product of any two singular
distributions which have their singularities at the same point on the unit
circle. If one considers integrating the resulting object and using for
this purpose the fundamental property of any of the two distributions
involved, one can see that the integral is not well defined in the context
of the definitions given here for the singular distributions. Although one
cannot rule out that some other definition can be found to include some
such cases, we certainly do not have one at this time.

We thus see that we are in fact unable to promote the whole space of
integrable real functions and singular distributions to an algebra.
However, there are some sub-spaces within which this can be done. Under
some circumstances one can solve the problem of defining within the
complex-analytic structure the product of two integrable real functions.
This cannot be done for the whole sub-space of integrable real functions,
because there is the possibility that the product of two integrable real
function will not be integrable. However, if we restrict the sub-space to
those integrable real functions which are also limited, then it can be
done. This is so because the product of two limited integrable real
functions is also a limited real function, and therefore integrable. In
this way, one can find the inner analytic function that corresponds to the
product, since according to the construction which was presented
in~\cite{CAoRFI}, this can be done for any integrable real function. The
resulting inner analytic function will not, however, be related in a
simple way to the inner analytic functions of the two factor functions.

One case in which the product can always be defined is that of an
integrable real function with a Dirac delta ``function'', so long as the
real function is well defined at the singular point of the delta
``function''. Given the nature of the delta ``function'', this is
equivalent to multiplying it by a mere real number, the value of the
integrable real function at the singular point of the delta ``function'',

\begin{equation}
  g(\theta)\delta(\theta-\theta_{1})
  =
  g(\theta_{1})\delta(\theta-\theta_{1}).
\end{equation}

\noindent
The corresponding inner analytic function is therefore given simply by
$g(\theta_{1})w_{\delta}(z,z_{1})$. Similar statements are true, of
course, for all the derivatives of the delta ``function''. Therefore, in
all such cases there is no difficulty in determining the inner analytic
function that corresponds to the product.

Note that this difficulty relates only to the definition of the product of
two real objects on the unit circle. As was observed before, for all the
singular distributions their definition by means of inner analytic
functions always provides the means to determine whether or not they can
be applied to a given real object, so long as it is represented by an
inner analytic function, and determines what results from that operation,
if it is possible at all.

\section{Conclusions and Outlook}\label{Sec06}

We have extended the close and deep relationship established in a previous
paper~\cite{CAoRFI}, between integrable real functions and complex
analytic functions in the unit disk centered at the origin of the complex
plane, to include singular distributions. This close relationship between,
on the one hand, real functions and related objects, and on the other
hand, complex analytic functions, allows one to use the powerful and
extremely well-known machinery of complex analysis to deal with the real
functions and related objects in a very robust way, even if these objects
are very far from being analytic. The concept of integral-differential
chains of proper inner analytic functions, which we introduced in that
previous paper, played a central role in the analysis presented.

One does not usually associate non-differentiable, discontinuous and
unbounded real functions with single analytic functions. Therefore, it may
come as a bit of a surprise that, as was established in~\cite{CAoRFI},
{\em all} integrable real functions are given by the real parts of certain
inner analytic functions on the open unit disk when one approaches the
unit circle. This surprise is now compounded by the fact that inner
analytic functions can represent singular distributions as well and, in
fact, can represent what may be understood as a complete set of such
singular objects.

There are many more inner analytic functions within the open unit disk
than those that were examined here and in~\cite{CAoRFI}, in relation to
integrable real functions and singular distributions. Therefore, it may be
possible to further generalize the relationship between real objects on
the unit circle and inner analytic functions. For example, we have
observed in this paper that there are inner analytic functions whose real
parts converge to non-integrable real functions on the unit circle. Simple
examples are the inner analytic functions given by

\begin{equation}
  \bar{w}_{\delta^{n\prime}}(z,z_{1})
  =
  -\ii 
  w_{\delta^{n\prime}}(z,z_{1}),
\end{equation}

\noindent
for $n\in\{0,1,2,3,\ldots,\infty\}$, that correspond to the non-integrable
real functions which are the Fourier-conjugate functions of the Dirac
delta ``function'' and its derivatives. This suggests that we consider the
question of how far this can be generalized, that is, of what is the
largest set of non-integrable real functions that can be represented by
inner analytic functions. This issue will be discussed in the fourth paper
of this series.

The singular distributions are integrable real objects associated to
non-integrable singularities of the corresponding inner analytic
functions, a fact which is made possible by the orientation of the
singularities with respect to the direction along the unit circle. This
suggests that the most general definition of such singular distributions
may be formulated in terms of the type and orientation of the
singularities present on the unit circle. In this case one would expect
that singular distributions would be associated to inner analytic
functions with hard singularities that are oriented in a particular way,
so that the integrals of their real parts can be defined via limits from
the open unit disk to the unit circle\footnote{Post-publication note:
  another way to state this is to say that the inner analytic functions
  associated to singular distributions must be {\em irregular} inner
  analytic functions, since they are not integrable around their singular
  points, which are hard singular points with strictly positive degrees of
  hardness, while the corresponding real objects are.}.

We believe that the results presented here establish a new perspective for
the representation and manipulation of singular distributions. It might
also constitute a simpler and more straightforward way to formulate and
develop the whole theory of Schwartz distributions within a compact
domain.

\section*{Acknowledgments}

The author would like to thank his friend and colleague Prof. Carlos
Eugênio Imbassay Carneiro, to whom he is deeply indebted for all his
interest and help, as well as his careful reading of the manuscript and
helpful criticism regarding this work.

\bibliography{allrefs_en}\bibliographystyle{ieeetr}

\end{document}

%% file: packages.tex
\usepackage[latin1]{inputenc}
\usepackage{epsfig}
\usepackage{amsfonts}
\usepackage{cite}
\usepackage{color}
\definecolor{lgrey}{gray}{0.99}

%% file: math-defs.tex
\newcommand{\ii}{\mbox{\boldmath$\imath$}}
\newcommand{\iii}{\mbox{\boldmath\scriptsize$\imath$}}
\newcommand{\Frac}[2]{\frac{\displaystyle #1}{\displaystyle #2}}
\newcommand{\FFrac}[2]{{\displaystyle\frac{\displaystyle #1}{\displaystyle #2}}}

\newcommand{\e}[1]{\,{\rm e}^{#1}}

\newcommand{\ldot}{\mbox{\Large$\cdot$}\!}

\newtheorem{definition}{Definition}

\newtheorem{theorem}{Theorem}

\newtheorem{proof}{Proof}[theorem]
\newcommand{\Colon}{{\hspace{-0.4em}\bf:}}